\documentclass[preprint,12pt]{elsarticle}

\usepackage{graphics}
\usepackage{graphicx}
\usepackage{epsfig}
\usepackage{amsmath}
\numberwithin{equation}{section}
\newtheorem{Def}{Definition}
\newtheorem{Thm}{Theorem}
\newtheorem{Rmk}{Remark}
\newtheorem{Cor}{Corollary}

\newtheorem{Pro}{Proposition}
\newtheorem{Ex}{Example}

\usepackage{amssymb}

\begin{document}

\begin{frontmatter}



\title{Ollivier's Ricci curvature, local clustering and \\curvature dimension inequalities on graphs}


\author[1,4,2]{J\"urgen Jost\fnref{fn1}}
\ead{jost@mis.mpg.de}
\author[1,3]{Shiping Liu\corref{cor1}}
\ead{shiping@mis.mpg.de}
\address[1]{Max Planck Institute for Mathematics in the Sciences,
Leipzig 04103, Germany}
\address[4]{Department of Mathematics and Computer Sciences, University of Leipzig, Leipzig 04103, Germany}
\address[2]{Santa Fe Institute for the Sciences of Complexity, Santa Fe,
NM 87501, USA}
\address[3]{Academy of Mathematics and Systems
Science, Chinese Academy of Sciences, Beijing 100190, China.}
\cortext[cor1]{Corresponding author.}
\fntext[fn1]{The research leading to these results has received funding from the
European Research Council under the European Union's Seventh
Framework Programme (FP7/2007-2013) / ERC grant agreement
n$^\circ$ 267087.}

\begin{abstract}
In this paper, we explore the relationship between one of the most
elementary and important properties of graphs, the presence and
relative frequency of triangles, and a combinatorial notion of Ricci
curvature. We employ
 a
  definition of generalized Ricci curvature proposed by Ollivier in
  a general framework of Markov processes and metric spaces and
  applied in graph theory by Lin-Yau. In analogy with curvature
  notions in Riemannian geometry, we interpret this Ricci curvature
  as a control on the amount of overlap between neighborhoods of two
  neighboring vertices. It is therefore naturally related to the
  presence of triangles containing  those vertices, or more precisely,
  the local clustering coefficient, that is,
  the relative proportion of connected neighbors among all the
  neighbors of a vertex.
  This suggests to derive lower Ricci curvature
  bounds on graphs in terms of such local clustering coefficients. We also study curvature dimension
  inequalities on graphs, building upon previous work of several
  authors.
\end{abstract}


\end{frontmatter}


\section{Introduction}
When one studies empirical graphs, one of the most obvious and basic
properties to investigate is the presence and number of triangles,
that is, connected triples of vertices. In bipartite graphs, for
instance, there are no triangles, whereas in a complete graph, every
triple of vertices constitutes a triangle. A basic observation then is
that when two neighboring vertices are contained in a triangle, their
neighborhoods of radius 1 (let's assign to every edge the length 1 for
the discussion in this introduction) share the third vertex of the
triangle. That is, the more triangles those two neighboring vertices
are contained in, the larger the overlap of their neighborhoods. This
suggests an analogy with the notion of Ricci curvature in Riemannian
geometry where a lower bound on the Ricci curvature also controls the
amount of overlaps of distance balls from below. This is what we are
going to explore in a quantitative manner in this paper.

In fact, Ricci curvature is a fundamental concept in Riemannian
geometry, see e.g. \cite{Jost}. It is
a quantity computed from second derivatives of the metric tensor. It
controls how fast geodesics starting at the same point diverge on
average. Equivalently, it controls how fast the volume of distance
balls grows as a function of the radius. As already indicated, it also controls the amount
of overlap of two distance balls in terms of their radii and the
distance between their centers. In fact, such lower bounds follow from
a lower bound on the Ricci curvature. It was then natural to look for
generalizations of such phenomena on metric spaces more general than
Riemannian manifolds. That is, the question to find substitutes for the
lower bounds on the above mentioned second derivative combinations of
the metric tensor that yield the same geometric control on a general
metric space. By now, there exist several insightful definitions of  synthetic Ricci curvature on general metric measure spaces, see
Sturm \cite{Sturm1, Sturm2}, Lott-Villani \cite{LV}, Ohta \cite{Ohta}, Ollivier \cite{Oll}
etc.

As indicated, in this paper, we want to explore the implications of such ideas in
graph theory. The geometric idea is that a lower Ricci curvature bound
prevents geodesics from diverging too fast and balls from growing too
fast in volume. On a graph, the analogue of geodesics starting in
different directions, but eventually approaching each other again,
would be a triangle. Therefore, it is natural that the Ricci curvature
on a graph should be related to the relative abundance of
triangles. The latter is captured by the local clustering coefficient
introduced by Watts-Strogatz \cite{WS}. Thus,  the intuition of
Ricci curvature on a graph should play with the relative frequency of
triangles a vertex shares with its neighbors. In fact, more precisely,
since the local clustering coefficient averages over the neighbors of
a vertex, this should really related to some notion of scalar
curvature, as an average of Ricci curvatures in different directions,
that is, for different neighbors of a given vertex.

Among the several definitions of
generalized Ricci curvature in the literature mentioned above, the one of
Ollivier works particularly well on
discrete spaces like graphs. It is formulated in terms of the
transportation distance between local measures:
\begin{equation}
\kappa (x, y):=1-W_1(m_x, m_y), \label{0}
\end{equation}
where $x,y$ are  vertices in our graph that are neighbors (written as
$x\sim y$) and the measure $m_x=\frac{1}{d_x}$, where $d_x$ is the
degree of $x$, puts equal weight on all neighbors. $W_1(m_x,m_y)$ is
the transportation distance between the two measures $m_x$ and $m_y$
(defined more precisely below).
When two balls
strongly overlap, as is the case in Riemannian geometry when the Ricci curvature has a large
lower bound, then it is easier to transport the mass of one to the
other. Analogously, in the graph case, when the two vertices share
many triangles, then the transportation distance should be smaller,
and the curvature therefore correspondingly larger. This is the idea of Ollivier's definition as we see it and
explore in this paper. We shall obtain both upper and lower bounds for
Ollivier's Ricci curvature on graphs in Section \ref{84}, which are
optimal on many graphs.

Let us now formulate our main result (recalled and proved below as
Theorem \ref{8}).
\begin{Thm} \label{t0}
 On a locally finite graph, we put for any pair of
 neighboring vertices $x, y$,
$$\sharp(x,y):=\text
{number of triangles which include}\,\,\, x, y \,\,\,\text{as vertices}=\sum_{x_1, x_1\sim x, x_1\sim y}1.$$
We then have
\begin{equation}\label{0a}
 \kappa (x,y)\geq-\left(1-\frac{1}{d_x}-\frac{1}{d_y}-\frac{\sharp(x,y)}{d_x\wedge d_y}\right)_+-\left(1-\frac{1}{d_x}-\frac{1}{d_y}-\frac{\sharp(x,y)}{d_x\vee d_y}\right)_++\frac{\sharp(x,y)}{d_x\vee d_y}.
\end{equation}
where  $s_+:=\max(s,0), s\vee t:=\max(s,t), s\wedge t:=\min(s,t)$.
\end{Thm}
This equality is sharp for instance for a complete graph of $n$
vertices where the left and the right hand side both equal to
$\frac{n-2}{n-1}$.

The local clustering coefficient introduced by Watts-Strogatz
\cite{WS} is
\begin{equation}
c(x):=\frac{\text{number of edges between neighbors of}\,\, x}{\text{number of possible existing edges between neighbors of}\,\, x},
\end{equation}
which measures the extent to which neighbors of $x$ are directly
connected, i.e.,
\begin{equation}\label{0c}
c(x)=\frac{1}{d_x(d_x-1)}\sum_{y, y\sim x}\sharp(x, y).
\end{equation}
Thus, this local clustering coefficient is an average over the
$\sharp(x,y)$ for the neighbors of $x$. Thus, we might also introduce
some kind of scalar curvature (suggested in Problem Q in Ollivier \cite{Oll2}) as
\begin{equation}\label{0d}
\kappa(x):= \frac{1}{d_x}\sum_{y, y\sim x}\kappa(x,y).
\end{equation}
For illustration, let us consider the case where our graph is $d$-regular, that is,
$d_z=d$ for all vertices $z$.
When $1\ge \frac{2}{d}+\frac{\sharp(x,y)}{d}$ for all $y\sim x$, we would then get
\begin{equation}
  \label{0e}
  \kappa(x)\ge -2 +\frac{4}{d} + \frac{3(d-1)}{d}c(x).
\end{equation}
This example nicely illustrates the relation between Ollivier's
curvature and the Watts-Strogatz clustering coefficient.

Without the triangle terms $\sharp(x,y)$ (which is the crucial term
for our purposes), Theorem \ref{t0} is due to
Lin-Yau \cite{LinYau,Lin}, and we take their proof as our starting
point. Lin-Yau also obtain analogues of Bochner type inequalities and
eigenvalue estimates as known from Riemannian geometry.

In Riemannian geometry, the Bochner formula encodes deep analytic
properties of Ricci curvature. It is a key ingredient in proving many
results, e.g. the spectral gap of the Laplace-Beltrami operator. A
lower bound of the Ricci curvature  implies a curvature dimension
inequality involving  the Laplace-Beltrami operator through the Bochner
formula. In an important work, Bakry and \'Emery \cite{BaEm, BaEm2}
generalize this inequality to generators of Markov semigroups, which
works on measure spaces. Their inequality contains plentiful
information and implies a lot of functional inequalities including
spectral gap inequalities, Sobolev inequalities, and logarithmic
Sobolev inequalities and many celebrated geometric theorems (see
\cite{Bak} and the references therein). Lin-Yau \cite{LinYau} study such inequalities on locally finite graphs.

In the present paper, we also want  to find relations on locally
finite graphs  between Ollivier's Ricci curvature and Bakry-\'Emery's
curvature dimension inequalities, which represent the geometric and
analytic aspects of graphs respectively. Again, this is inspired  by
Riemannian geometry where one may attach a Brownian motion with a drift
to a Riemannian metric \cite{Oll}. We also mention that the definitions
given by Sturm and Lott-Villani are also consistent with that of
Bakry-\'Emery \cite {Sturm1, Sturm2, LV}. So exploring the relations
on nonsmooth spaces may provide a good point of view to connect
Ollivier's definition to Sturm and Lott-Villani's (in this aspect, see also Ollivier-Villani \cite{OllVill}). In Section \ref{91},
we use the local clustering coefficient again to establish more
precise curvature dimension inequalities than those of Lin-Yau \cite{LinYau}. And with this in hand, we prove curvature dimension inequalities under the condition that Ollivier's Ricci curvature of the graph is positive.

Further
analytical results following  from curvature dimension inequalities on
finite graphs have been described in \cite{Lin}, and Lin-Lu-Yau
\cite{LLY} study a modified definition of Ollivier's Ricci curvature
on graphs. Recently, Paeng \cite{Paeng} studied upper bounds for the  diameter and volume of finite simple graphs in terms of Ollivier's Ricci curvature. For other works of synthetic Ricci curvatures on discrete spaces, see Dodziuk-Karp \cite{DK}, Chung-Yau \cite{CY}, Bonciocat-Sturm \cite{BS}, and on cell complexes see Forman \cite{Forman}, Stone \cite{Stone} etc.

We point out that, as in Riemannian geometry, both Ollivier's Ricci
curvature and Bakry-\'Emery's curvature dimension inequality can
give lower bound estimates of the first eigenvalue $\lambda_1$ for the
Laplace operator (see Ollivier \cite{Oll}, Bakry
\cite{Bak}). Therefore our results in fact relate $\lambda_1$ to the
Watts-Strogatz local clustering coefficient, or the number of cycles
with length $3$. In \cite{DS}, Diaconis and Stroock obtain several
geometric bounds for eigenvalues of graphs, one of which is related to
the number of odd length cycles. For more geometric quantities and
methods concerning eigenvalue estimates in the study of Markov chains,
see \cite{DSC, D1, D2} and the references therein. We  further
explore the interaction between Ollivier's Ricci curvature and
eigenvalues estimates in joint work with Frank Bauer, see \cite{BJL}.

In this paper, $G=(V, E)$ will denote an undirected connected simple
graph without loops, where $V$ is the set of vertices and $E$ is the
set of edges. $V$ could be an infinite set. But we require that $G$ is
locally finite, i.e., for every $x\in V$, the number of edges
connected to $x$ is finite. For simplicity and in order to see more geometry, we mainly work
on unweighted graphs. But we will also derive similar results on
weighted graphs. In that case, we denote by $w_{xy}$ the weight
associated to $x, y\in V$, where $x\sim y$ (we may simply put
$w_{xy}=0$ if $x$ and $y$ are not neighbors, to simplify the
notation). The unweighted case corresponds to $w_{xy}=1$ whenever $ x\sim y$. The degree of $x\in V$ is $d_x=\sum_{y, y\sim x}w_{xy}$.

\section{Ollivier's Ricci curvature and Bakry-\'Emery's calculus}
In this section, we present some basic facts about Ollivier's Ricci
curvature and Bakry-Emery's $\Gamma_2$ calculus, in particular on graphs.

\subsection{Ollivier's Ricci curvature}

Ollivier's Ricci curvature works on a general metric space $(X, d)$, on which we attach to each point $x\in X$ a probability measure $m_x(\cdot)$. We denote this structure by $(X, d, m)$.

For a locally finite unweighted graph $G=(V, E)$, we define the metric $d$ as follows. For neighbors $x, y$, $d(x, y)=1$. For general distinct vertices $x, y$,  $d(x, y)$ is the length of the shortest path connecting $x$ and $y$, i.e. the number of edges of the path. We attach to each vertices $x\in V$ a probability measure
\begin{equation}\label{85}
 m_x(y)=\left\{
          \begin{array}{ll}
            \frac{1}{d_x}, & \hbox{if $y \sim x$;} \\
            0, & \hbox{otherwise.}
          \end{array}
        \right.
\end{equation}
An intuitive illustration of this is a random walker that sits at $x$
and then chooses amongst the neighbors of $x$ with equal probability
$\frac{1}{d_x}$.
\begin{Def}[Ollivier]
On $(X, d, m)$, for any two distinct points $x, y\in X$, the (Ollivier-) Ricci curvature of $(X, d, m)$ along $(xy)$ is defined as
\begin{equation}
\kappa (x, y):=1-\frac{W_1(m_x, m_y)}{d(x, y)}.\label{1}
\end{equation}
\end{Def}
Here, $W_1(m_x, m_y)$ is the optimal transportation distance
 between the two
probability measures $m_x$ and $m_y$, defined as follows (cf. Villani \cite{V1, V2}, Evans \cite{Evans}).
\begin{Def}
For two probability measures $\mu_1$, $\mu_2$ on a metric space $(X, d)$, the transportation
distance between them is defined as
\begin{equation}
W_1(\mu_1, \mu_2):= \inf_{\xi\in\prod(\mu_1, \mu_2)} \int_{X\times X}d(x,y)d\xi(x,y),\label{2}
\end{equation}
where $\prod(\mu_1, \mu_2)$ is the set of probability measures on $X\times X$ projecting to $\mu_1$ and $\mu_2$.
\end{Def}
In other words, $\xi$ satisfies
\begin{equation*}
\xi(A\times X)=\mu_1(A),\,\,\,\xi(X\times B)=\mu_2(B),\,\,\forall A, B\subset X.
\end{equation*}
\begin{Rmk}
Intuitively, this distance measures the optimal cost to move one pile
of sand to another one with the same mass. For  case of a graph $G=(G,
d, m)$, the supports of $m_x$ and $m_y$ are  finite discrete sets, and
thus, $\xi$ is just a matrix with terms $\xi (x', y')$ representing the mass moving from $x'\in \text{support of}\,\, m_x$ to $y'\in \text{support of} \,\,m_y$.  That is, in this case,
$$W_1(m_x, m_y)=\inf_{\xi}\sum_{x', x'\sim x}\sum_{y', y'\sim y}d(x', y')\xi(x', y'),$$
where the infimum is taken over all matrices $\xi$ which satisfy
$$\sum_{x', x'\sim x}\xi(x', y')=\frac{w_{yy'}}{d_y},\,\,\,\sum_{y', y'\sim y}\xi(x', y')=\frac{w_{xx'}}{d_x}.$$
We also call $\xi$ a transfer plan. If we can find a particular transfer plan, we then get an upper bound for $W_1$ and therefore a lower bound for $\kappa$.
\end{Rmk}

A very important property of transportation distance is the Kantorovich duality (see, e.g. Theorem 1.14 in Villani \cite{V1}). We state it here in our particular graph setting

\begin{Pro}[Kantorovich duality]
$$W_1(m_x, m_y)=\sup_{f, 1-Lip}\left[\sum_{z, z\sim x} f(z) m_x(z)-\sum_{z, z\sim y} f(z)dm_y(z)\right],$$
where the supremum is taken over all functions on $G$ that satisfy
$$|f(x)-f(y)|\leq d(x, y),$$
for any $x, y\in V$, $x\neq y$.
\end{Pro}

From this property, a good choice of a $1$-Lipschitz function $f$ will
yield a lower bound for $W_1$ and therefore an upper bound for $\kappa$.
\begin{Rmk}
We list some basic first observations about this curvature concept (see Ollivier \cite{Oll}):
\begin{itemize}
  \item $\kappa(x, y)\leq 1$.
  \item Rewriting  (\ref{1}) gives $W_1(m_x, m_y)=d(x,y)(1-\kappa (x,
    y))$, which is analogous to the expansion in the Riemannian case.
  \item A lower bound $\kappa (x, y)\geq k$ for any $x, y\in X$ implies
  \begin{equation}
  W_1(m_x, m_y)\leq (1-k)d(x,y),\label{3}
  \end{equation}
  which can be seen as some kind of Lipschitz continuity of  measures.
\end{itemize}
\end{Rmk}

\subsection{Bakry-\'Emery's curvature-dimension inequality}
\subsubsection{Laplace operator}
We will study the following operator which is an analogue of the Laplace-Beltrami operator in Riemannian geometry.
\begin{Def}
The Laplace operator on $(X,d,m)$ is defined as follows
\begin{equation}\label{95}
  \Delta f(x)=\int_X f(y)dm_x(y)-f(x), \,\,\text{for functions}\,\,\, f: X\longrightarrow R.
\end{equation}
\end{Def}

For our choice of $\{m_x(\cdot)\}$, this is the graph Laplacian
studied by many authors, see e.g. \cite{Chung}, \cite{BanJ}, \cite{BJ}, \cite{DK}, \cite{LinYau}.

\subsubsection{Bochner formula and curvature-dimension inequality}
In the Riemannian case, many analytical consequences of a lower bound
of the Ricci curvature are obtained through the well-known Bochner formula,
$$\frac{1}{2}\Delta(|\nabla f|^2)=|Hess\,\, f|^2+\langle\nabla(\Delta f),\nabla f\rangle+Ric(\nabla f,\nabla f).$$
Analytically, $|Hess\,\, f|^2$ is difficult to be defined on a nonsmooth space. But using Schwarz's inequality, we have
$$|Hess\,\, f|^2\geq\frac{(\Delta f)^2}{m},$$ where $m$ is the dimension constant. So we can use
\begin{equation}\label{83}
\frac{1}{2}\Delta(|\nabla f|^2)\geq\frac{(\Delta f)^2}{m}+\langle\nabla(\Delta f),\nabla f\rangle+K|\nabla f|^2
\end{equation}
to characterize $Ric \geq K$.

Bakry-\'Emery \cite{Bak, BaEm, BaEm2} take this inequality as the starting
point and directly use the
operators to define curvature bounds. Starting from an operator $\Delta$, they define iteratively,
\begin{align*}
\Gamma_0(f, g)&=fg,\\\Gamma(f, g)&=\frac{1}{2}\{\Delta\Gamma_0(f, g)-\Gamma_0(f, \Delta g)-\Gamma_0(\Delta f, g)\},\\
\Gamma_2(f, g)&=\frac{1}{2}\{\Delta\Gamma(f, g)-\Gamma(f, \Delta g)-\Gamma(\Delta f, g)\}.
\end{align*}
In fact, $\Gamma(f, f)$ is an analogue of $|\nabla f|^2$, and $\Gamma_2(f, f)$ is an analogue of $\frac{1}{2}\Delta|\nabla f|^2-\langle\nabla(\Delta f), \nabla f\rangle$ in (\ref{83}).
\begin{Def}
We say an operator $\Delta$ satisfies a curvature-dimension inequality $CD(m, K) $ if for all functions $f$ in the domain of the operator
\begin{equation}
\Gamma_2(f, f)(x)\geq \frac{1}{m}(\Delta f(x))^2+K(x)\Gamma(f, f)(x), \,\,\forall x\in X,
\end{equation}
where $m\in[1, +\infty]$ is the dimension parameter, $K(x)$ is the curvature function.
\end{Def}

As studied in Lin-Yau \cite{LinYau}, applying this construction to
the operator (\ref{95}) gives
\begin{align}\label{87}
\Gamma(f, f)(x)&=\frac{1}{2}\int_X(f(y)-f(x))^2dm_x(y),
\end{align}
In fact generally
\begin{equation*}
\Gamma(f, g)(x)=\frac{1}{2}\int_X(f(y)-f(x))(g(y)-g(x))dm_x(y).
\end{equation*}

For the sake of convenience, we will denote
$$Hf(x):=\frac{1}{4}\int_X\int_X(f(x)-2f(y)+f(z))^2dm_y(z)dm_x(y).$$

By the calculation in Lin-Yau \cite{LinYau} we get
\begin{align*}
\Delta\Gamma(f, f)(x)&=2Hf(x)-\int_X\int_X(f(x)-2f(y)+f(z))(f(x)-f(y))dm_y(z)dm_x(y),\\
2\Gamma(f, \Delta f)(x)&=-(\Delta f(x))^2-\int_X\int_X(f(z)-f(y))(f(x)-f(y))dm_y(z)dm_x(y).
\end{align*}
and then,
\begin{equation}\label{86}
\Gamma_2(f, f)=Hf(x)-\Gamma(f, f)(x)+\frac{1}{2}(\Delta f(x))^2.
\end{equation}

\section{Ollivier's Ricci curvature and triangles}\label{84}

In this section, we mainly prove lower bounds for Ollivier's Ricci
curvauture on locally finite graphs. In particular we shall explore
the implication between lower bounds of the curvature and the number
of triangles including neighboring vertices; the latter is encoded in
the local clustering coefficient. We remark that we only need to bound
$\kappa(x, y)$ from below for neighboring $x, y$, since by the
triangle inequality of $W_1$, this will also be a lower bound for
$\kappa(x,y)$ of any pair of $x, y$. (See Proposition 19 in Ollivier \cite{Oll}.)

\subsection{Unweighted graphs}
In this subsection, we only consider unweighted graphs.

In Lin-Yau \cite{LinYau}, they prove a lower bound of Ollivier's Ricci
curvature on locally finite graphs $G$. Here, for later purposes, we include the case where $G$ may have vertices of degree $1$ and get the following modified result.

\begin{Thm}\label{7}
On a locally finite graph $G=(V, E)$, we have for any pair of neighboring vertices $x, y$,
\begin{equation*}
\kappa (x,y)\geq-2\left(1-\frac{1}{d_x}-\frac{1}{d_y}\right)_+=\left\{
                                                                \begin{array}{ll}
                                                                  -2+\frac{2}{d_x}+\frac{2}{d_y}, & \hbox{if $d_x>1$ and $d_y>1$;} \\
                                                                  0, & \hbox{otherwise.}
                                                                \end{array}
                                                              \right.
\end{equation*}
\end{Thm}

\begin{Rmk}
Notice that if $d_x=1$, then we can calculate $\kappa(x, y)=0$
exactly. So, even though  in this case  $-2+\frac{2}{d_x}=0$, $\kappa(x, y)\geq \frac{2}{d_y}$ doesn't hold.
\end{Rmk}

For completeness, we state the proof of Theorem \ref{7} here. It is essentially the one in Lin-Yau \cite{LinYau} with a small modification.

\verb"Proof of Theorem" \ref{7}: Since $d(x,y)=1$ for $x\sim y$, we have
\begin{equation}
\kappa(x,y)=1-W_1(m_x, m_y).\label{21}
\end{equation}
Using Kantorovich duality, we get
\begin{align}
W_1(m_x, m_y)=&\sup_{f, 1-Lip}\left(\frac{1}{d_x}\sum_{z, z\sim x}f(z)-\frac{1}{d_y}\sum_{z', z'\sim y}f(z')\right)\nonumber
\\=&\sup_{f, 1-Lip}\Bigg(\frac{1}{d_x}\sum_{z, z\sim x, z \neq y}(f(z)-f(x))-\frac{1}{d_y}\sum_{z', z'\sim y, z'\neq x}(f(z')-f(y))\nonumber\\&+\frac{1}{d_x}(f(y)-f(x))-\frac{1}{d_y}(f(x)-f(y))+(f(x)-f(y))\Bigg)
\nonumber\\\leq &\frac{d_x-1}{d_x}+\frac{d_y-1}{d_y}+\left|1-\frac{1}{d_x}-\frac{1}{d_y}\right|\nonumber
\\=&2-\frac{1}{d_x}-\frac{1}{d_y}+\left|1-\frac{1}{d_x}-\frac{1}{d_y}\right|.\nonumber
\\=&1+2\left(1-\frac{1}{d_x}-\frac{1}{d_y}\right)_+.\label{80}
\end{align}
Inserting the above estimate into (\ref{21}) gives
\begin{align*}
\kappa(x, y)&\geq-2\left(1-\frac{1}{d_x}-\frac{1}{d_y}\right)_+.
\end{align*}$\hfill\Box$

Note that trees attain this lower bound. This coincides with the geometric intuition of curvature. Since trees have the fastest volume growth rate, it is plausible that they have the smallest curvature.

\begin{Pro}
 We consider a tree $T=(V, E)$. Then for any neighboring $x, y$, we have
\begin{equation}
\kappa (x,y)=-2\left(1-\frac{1}{d_x}-\frac{1}{d_y}\right)_+.
\end{equation}
\end{Pro}

\verb"Proof:" In fact with Theorem \ref{7} in hand, we only need to
prove that $1+2\left(1-\frac{1}{d_x}-\frac{1}{d_y}\right)_+$ is also a lower bound of $W_1$. If one of $x, y$ is a vertex of degree $1$, say $d_x=1$, it is obvious that $W_1(m_x, m_y)=1$. So we only need to deal with the case $1-\frac{1}{d_x}-\frac{1}{d_y}\geq 0$.

We can find a 1-Lipschitz function $f$ on a tree as follows.
\begin{equation}
 f(z)=\left\{
        \begin{array}{ll}
          0, & \hbox{if $z\sim y, z\neq x$;} \\
          1, & \hbox{if $z=y$;} \\
          2, & \hbox{if $z=x$;} \\
          3, & \hbox{if $z\sim x, z\neq x$.}
        \end{array}
      \right.
\end{equation}
Since on a tree, the path joining two vertices are unique, there is no further path between neighbors of $x$ and $y$. So this can be easily extended to a $1$-Lipschitz function on the whole graph. Then by Kantorovich duality, we have
\begin{align}
 W_1(m_x, m_y)&\geq \frac{1}{d_x}(3(d_x-1)+1)-\frac{1}{d_y}\cdot 2\nonumber\\
&=3-\frac{2}{d_x}-\frac{2}{d_y}.\label{10}
\end{align}

This completes the proof. $\hfill\Box$

In order to make clear the geometric meaning of the term
$\left(1-\frac{1}{d_x}-\frac{1}{d_y}\right)_+$, and  also to prepare the idea used in the next theorem, we give another method to get the upper bound of $W_1$. That works through a particular transfer plan.
If $$1-\frac{1}{d_x}-\frac{1}{d_y}\geq 0, \,\,\text{or}\,\,1-\frac{1}{d_y}\geq\frac{1}{d_x},$$
then for $m_y$, the mass at all $z$ such that $z\sim y$, $z\neq x$ is larger than that of $m_x$ at $y$. So we can move the mass $\frac{1}{d_x}$ at $y$ to $z$, $z\sim y$, $z\neq x$ for distance $1$. Symmetrically, we can move a mass of $\frac{1}{d_y}$ at the vertices $z$ which satisfy $z\sim x$, $z\neq y$ to $x$ for distance $1$. The remaining mass of $\left(1-\frac{1}{d_x}-\frac{1}{d_y}\right)$ needs to be moved for distance 3. This gives
\begin{align}
 W_1(m_x, m_y)&\leq \left(\frac{1}{d_x}+\frac{1}{d_y}\right)\times 1+\left(1-\frac{1}{d_x}-\frac{1}{d_y}\right)\times 3\nonumber\\&=3-\frac{2}{d_x}-\frac{2}{d_y}.\label{9}
\end{align}

If $$1-\frac{1}{d_x}-\frac{1}{d_y}\leq 0,$$ we only need to move the mass of $m_x$ for distance $1$ to the support of $m_y$. So we have in this case $W_1(m_x, m_y)=1$. This gives the same upper bound as in (\ref{80}).

From the view of transfer plans, the existence of triangles including neighboring vertices would save a lot of transport costs and therefore affect the curvature heavily. We denote for $x\sim y$,
$$\sharp(x,y):=\text{number of triangles which include}\,\,\, x, y \,\,\,\text{as vertices}=\sum_{x_1, x_1\sim x, x_1\sim y}1.$$

\begin{Rmk}
This quantity $\sharp(x, y)$ is related to  the local clustering coefficient introduced by Watts-Strogatz \cite{WS},
\begin{equation*}
c(x):=\frac{\text{number of edges between neighbors of}\,\, x}{\text{number of possible existing edges between neighbors of}\,\, x},
\end{equation*}
which measures the extent to which neighbors of $x$ are directly connected. In fact, we have the relation
\begin{equation}\label{81}
c(x)=\frac{1}{d_x(d_x-1)}\sum_{y, y\sim x}\sharp(x, y).
\end{equation}
\end{Rmk}

We will explore the relation between the curvature $\kappa(x, y)$ and
the number of triangles $\sharp(x, y)$. A critical observation is that $\kappa(x,y)$ is symmetric w.r.t. $x$ and $y$. So we try to express the curvature through symmetric quantities
$$d_x\wedge d_y:=\min \{d_x, d_y\},\,\,\, d_x\vee d_y:=\max\{d_x, d_y\}.$$

\begin{Thm}\label{8}
 On a locally finite graph $G=(V, E)$, we have for any pair of neighboring vertices $x, y$,
\begin{equation*}
 \kappa (x,y)\geq-\left(1-\frac{1}{d_x}-\frac{1}{d_y}-\frac{\sharp(x,y)}{d_x\wedge d_y}\right)_+-\left(1-\frac{1}{d_x}-\frac{1}{d_y}-\frac{\sharp(x,y)}{d_x\vee d_y}\right)_++\frac{\sharp(x,y)}{d_x\vee d_y}.
\end{equation*}
Moreover, this inequality is sharp for certain graphs.
\end{Thm}

\begin{Rmk}
If $\sharp(x, y)=0$, then this lower bound reduces to the one in Theorem \ref{7}.
\end{Rmk}
\begin{Ex}\label{71}
On a complete graph ${\cal K}_n$ $(n\geq 2)$ with $n$ vertices, $\sharp(x, y)=n-2$ for any $x, y$. So Theorem \ref{8} implies
$$\kappa(x, y)\geq \frac{n-2}{n-1}.$$
In fact, we can easily check that the above inequality is an equality. Also notice that on those graphs, the local clustering coefficient $c(x)=1$ attains the largest value.
\end{Ex}

Before carrying out the proof of Theorem \ref{8}, we  fix some notations. The vertices $z$ that are adjacent to
$x$ or $y$, where $x\sim y$, are divided into three classes.
\begin{itemize}
  \item common neighbors of $x, y$: $z\sim x$ and $z\sim y$:
  \item $x$'s own neighbors: $z\sim x, z\not\sim y, z\neq y$;
  \item $y$'s own neighbors: $z\sim y, z\not\sim x, z\neq x$.
\end{itemize}
\verb"Proof of Theorem" \ref{8}: We suppose w.l.o.g.,
$$d_x=d_x\vee d_y,\,\,\,d_y=d_x\wedge d_y.$$
In principle, our transfer plan moving $m_x$ to $m_y$ should be as follows.

\begin{enumerate}
\item Move the mass of $\frac{1}{d_x}$ from $y$ to $y$'s own neighbors;
\item Move a mass of $\frac{1}{d_y}$ from $x$'s own neighbors to $x$;
\item Fill gaps using the mass at $x$'s own neighbors. Filling the gaps at common neighbors costs $2$ and the one at $y$'s own neighbors costs $3$.
\end{enumerate}
A critical point will be whether $(1)$ and $(2)$ can be realized or not. It is easy to see that we can realize step $(1)$ if and only if
\begin{equation}
1-\frac{1}{d_y}-\frac{\sharp(x, y)}{d_y}\geq \frac{1}{d_x},\,\,\,\text{or}\,\,\,\,A:=1-\frac{1}{d_x}-\frac{1}{d_y}-\frac{\sharp(x, y)}{d_x\wedge d_y}\geq 0.
\end{equation}
That is, after taking off the mass at $x$ and common neighbors, $m_y$ still has at least a mass of $\frac{1}{d_x}$. Step $(2)$ can be realized if and only if
\begin{equation}
1-\frac{1}{d_x}-\frac{\sharp(x, y)}{d_x}\geq \frac{1}{d_y},\,\,\,\text{or}\,\,\,\,B:=1-\frac{1}{d_x}-\frac{1}{d_y}-\frac{\sharp(x, y)}{d_x\vee d_y}\geq 0.
\end{equation}
That is, after taking off the mass at $y$ and common neighbors, $m_x$ still has enough mass to fill $\frac{1}{d_y}$. Obviously, $A\leq B$.

We will divide the discussion into $3$ cases according to whether the first two steps can be realized or not.

\begin{itemize}
\item $0\leq A\leq B$.
This means we can adopt the above transfer plan. By definition of $W_1(m_x, m_y)$, we get
\begin{align*}
W_1(m_x, m_y)\leq & \frac{1}{d_x}\times 1+\frac{1}{d_y}\times 1+\left(\frac{1}{d_y}-\frac{1}{d_x}\right)\times \sharp(x, y)\times 2\\&+\left[1-\frac{1}{d_x}-\frac{1}{d_y}-\left(\frac{1}{d_y}-\frac{1}{d_x}\right)\times \sharp(x, y)-\frac{1}{d_x}\sharp(x, y)\right]\times 3\\=&3-\frac{2}{d_x}-\frac{2}{d_y}-\frac{\sharp(x, y)}{d_y}-\frac{2\sharp(x, y)}{d_x}.
\end{align*}
Or in a symmetric way,
\begin{equation}
W_1(m_x, m_y)\leq 3-\frac{2}{d_x\vee d_y}-\frac{2}{d_x\wedge d_y}-\frac{\sharp(x, y)}{d_x\wedge d_y}-\frac{2\sharp(x, y)}{d_x\vee d_y}.\label{99}
\end{equation}
Moreover, in this case the following function $f$ (as shown in Figure $1$) can be extended as a $1$-Lipschitz function,

\setlength{\unitlength}{1.5cm}
\begin{center}
\begin{picture}(4,2.25)\thicklines
\put(1,1){\line(1,0){2.0}}
\put(1,1){\line(1,1){1.0}}
\put(3,1){\line(-1,1){1.0}}
\put(1,1){\line(-1,1){0.8}}
\put(1,1){\line(-1,-1){0.7}}
\put(1,1){\line(-1,-2){0.5}}
\put(3,1){\line(1,1){0.6}}
\put(3,1){\line(1,-1){0.6}}
\put(1.1,0.9){\makebox(0,0){x}}
\put(2.9,0.85){\makebox(0,0){y}}
\put(1,1.2){\makebox(0,0){2}}
\put(0.2,1.95){\makebox(0,0){3}}
\put(0.3, 0.15){\makebox(0,0){3}}
\put(0.7,0){\makebox(0,0){3}}
\put(2,2.15){\makebox(0,0){1}}
\put(3,1.2){\makebox(0,0){1}}
\put(3.6,1.75){\makebox(0,0){0}}
\put(3.6,0.25){\makebox(0,0){0}}
\end{picture}

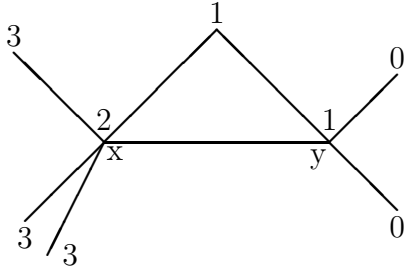
Figure 1. Mass moved from vertices with larger value

to those with smaller ones.
\end{center}

\begin{equation*}
f(z)=\left\{\begin{array}{ll}
            0, & \hbox{at $y$'s own neighbors;}\\
            1, & \hbox{at $y$ or common neighbors;}\\
            2, & \hbox{at $x$;}\\
            3, & \hbox{at $x$'s own neighbor,}
            \end{array}\right.
\end{equation*}
(that is, if there are no paths of length $1$ between common neighbors and $x$'s own neighbors,  nor paths of length $1$ or $2$ between $x$'s own neighbors and $y$'s own ones,) we have by Kantorovich duality,
\begin{align*}
W_1(m_x, m_y)\geq& \frac{1}{d_x}[f(y)+3(d_x-1-\sharp(x, y))+\sharp(x, y)]-\frac{1}{d_y}(f(x)+\sharp(x, y))\\
=&3-\frac{2}{d_x}-\frac{2}{d_y}-\frac{\sharp(x, y)}{d_y}-\frac{2\sharp(x, y)}{d_x}.
\end{align*}
That is, in this case, (\ref{9}) should be an equality.
In conclusion, $$\kappa(x, y)\geq -2+\frac{2}{d_x}+\frac{2}{d_y}+\frac{\sharp(x, y)}{d_x\wedge d_y}+\frac{2\sharp(x, y)}{d_x\vee d_y},$$
and the "$=$" can be attained.
\begin{Rmk}
$A\geq 0$ is equivalent to
\begin{equation*}
d_x\wedge d_y>1,\,\,\text{and}\,\,\,\,\sharp(x, y)\leq d_x\wedge d_y-1-\frac{d_x\wedge d_y}{d_x\vee d_y}.
\end{equation*}
Since $\sharp(x, y)\in \textbf{Z}$, we know that $d_x\wedge d_y\geq 2$ and $\sharp(x, y)\leq d_x\wedge d_y-2$. This means both $x$ and $y$ have at least one own neighbor.

If $A<0$, we get
$$d_x\wedge d_y-1-\frac{d_x\wedge d_y}{d_x\vee d_y}< \sharp(x, y)\leq d_x\wedge d_y-1.$$
I.e., $\sharp(x, y)=d_x\wedge d_y-1$. This means the vertex with smaller degree has no own neighbors.
\end{Rmk}
\item $A<0\leq B$.
In this case we cannot realize  step $(1)$ but step $(2)$ can be
realized. By the above remark, $A<0$ implies that  $y$ has no own neighbors. Our transfer plan should be step $(2)$ at first. Since $B\geq 0$ also implies
\begin{equation}
1-\frac{1}{d_y}-\frac{\sharp(x, y)}{d_x}\geq\frac{1}{d_x},
\end{equation}
so we can move the mass of $\frac{1}{d_x}$ at $y$ for distance $1$ to
common neighbors. Finally, we fill the gap at common neighbors for
distance $2$. In a formula,
\begin{align*}
W_1(m_x, m_y)\leq& \frac{1}{d_x}\times 1+\frac{1}{d_y}\times 1+\left(1-\frac{1}{d_x}-\frac{1}{d_y}-\frac{\sharp(x, y)}{d_x}\right)\times 2\\=& 2-\frac{1}{d_x}-\frac{1}{d_y}-\frac{2\sharp(x, y)}{d_x}.
\end{align*}
Or in a symmetric manner,
\begin{equation}
W_1(m_x, m_y)\leq 2-\frac{1}{d_x\vee d_y}-\frac{1}{d_x\wedge d_y}-\frac{2\sharp(x, y)}{d_x\vee d_y}.
\end{equation}
Moreover, in case the following function $f$ can be extended as a $1$-Lipschitz one,
\begin{equation*}
f(z)=\left\{\begin{array}{ll}
            0, & \hbox{at common neighbors;}\\
            1, & \hbox{at $x$ and $y$;}\\
            2, & \hbox{at $x$'s own neighbor,}
            \end{array}\right.
\end{equation*}
(that is, if there are no  paths of length $1$ between common neighbors and $x$'s own neighbors,) we have by Kantorovich duality,
\begin{align*}
W_1(m_x, m_y)\geq& \frac{1}{d_x}[f(y)+2(d_x-1-\sharp(x, y))]-\frac{1}{d_y}f(x)\\
=&2-\frac{1}{d_x}-\frac{1}{d_y}-\frac{2\sharp(x, y)}{d_x}.
\end{align*}
In conclusion, $$\kappa(x, y)\geq -1+\frac{1}{d_x}+\frac{1}{d_y}+\frac{2\sharp(x, y)}{d_x\vee d_y},$$
and the "$=$" can be attained.
\begin{Rmk}
Noting that if $\sharp(x, y)=d_x\wedge d_y-1$ then $B\geq 0$ is equivalent to
\begin{equation}\label{82}
d_x\vee d_y\geq \frac{d_x\wedge d_y}{d_x\wedge d_y-1}d_x\wedge d_y.
\end{equation}
In this case, one of $d_x$, $d_y$ has no own neighbors, and if  the
other one has sufficiently many own neighbors, $B\geq 0$ will be satisfied.
\end{Rmk}
\item $A\leq B<0$. In this case, neither step $(1)$ nor $(2)$ is
  applicable. Also,  $y$ has no own neighbor, and $B<0$ implies that
  we can move all the mass at $x$'s own neighbors to $x$ at first. And
  then we move the mass of $\frac{1}{d_x}$ at $y$ for distance $1$ to
  fill the gaps at $x$ and the common neighbors. In a formula,
\begin{align*}
W_1(m_x, m_y)\leq& \left(1-\frac{\sharp(x, y)}{d_x}\right)\times 1=1-\frac{\sharp(x, y)}{d_x}.
\end{align*}
Or in a symmetric way,
\begin{equation}
W_1(m_x, m_y)\leq 1-\frac{\sharp(x, y)}{d_x\vee d_y}.
\end{equation}
We can find a $1$-Lipschitz function
\begin{equation*}
f(z)=\left\{\begin{array}{ll}
            0, & \hbox{at $x$ and common neighbors;}\\
            1, & \hbox{at $y$ and $x$'s own neighbors,}
            \end{array}\right.
\end{equation*}
Then by Kantorovich duality,
\begin{align*}
W_1(m_x, m_y)\geq& \frac{1}{d_x}(f(y)+d_x-1-\sharp(x, y))-\frac{1}{d_y}\times 0\\
=&1-\frac{\sharp(x, y)}{d_x}.
\end{align*}
In this case $f$ can be extended to a $1$-Lipschitz function on the graph, so we get finally,
$$\kappa(x, y)=\frac{\sharp(x, y)}{d_x}.$$
\end{itemize}
Luckily, we can write those three cases in a uniform formula.$\hfill\Box$

\begin{Rmk}
From  extending $f$ to a $1$-Lipschitz function, we see that the paths
of length $1$ or $2$ between neighbors of $x$ and $y$ have an
important effect on the curvature. That is, in addition to triangles, quadrangles and pentagons are also related to Ollivier's Ricci curvature. But  polygons with more than 5 edges do not impact it.
\end{Rmk}

\begin{Rmk}
If we see the graph $G=(V, E)$ as a metric measure space $(G, d, m)$,
then the term $\sharp(x, y)/d_x\vee d_y$ is exactly $m_x\wedge
m_y(G):=m_x(G)-(m_x-m_y)_+(G)$, i.e. the intersection measure of $m_x$
and $m_y$. From a metric view, the vertices $x_1$ that satisfy
$x_1\sim x$, $x_1\sim y$ constitute the intersection of the  unit metric spheres $S_x(1)$ and $S_y(1)$.
\end{Rmk}

From Theorem \ref{8}, we can force the curvature $\kappa(x, y)$ to be positive by increasing the number $\sharp(x, y)$.

\begin{Thm}\label{89}
On a locally finite graph $G=(V, E)$, for any neighboring $x, y$, we have
\begin{equation}
 \kappa(x, y)\leq \frac{\sharp(x, y)}{d_x\vee d_y}.
\end{equation}
\end{Thm}
\verb"Proof:" Since except for the mass at common neighbors which we need not  move, the others have to be moved for a distance at least $1$, we have
$$W_1(m_x, m_y)\geq \left(1-\frac{\sharp(x, y)}{d_x\vee d_y}\right)\times 1.$$$\hfill\Box$

So if $\kappa(x, y)>0$, then $\sharp(x, y)$ is at least $1$. Moreover, if $\kappa(x, y)\geq k>0$, we have
\begin{equation}
 \sharp(x, y)\geq \lceil kd_x\vee d_y\rceil,
\end{equation}
where $\lceil a\rceil:=\min\{A\in \textbf{Z}|A\geq a\}$, for $a\in R$.

We will denote $D(x):=\max_{y, y\sim x}d_y$. By the relation (\ref{81}), we can get immediately
\begin{Cor}
The scalar curvature at $x$ can be controlled by the local clustering coefficient at $x$,
\begin{equation*}
\frac{d_x-1}{d_x}c(x)\geq\kappa(x)\geq-2+\frac{d_x-1}{d_x\vee D(x)}c(x).
\end{equation*}
\end{Cor}
\begin{Rmk}
In fact in some special cases, we can get more precise lower bounds
\begin{equation*}
\kappa(x)\geq\left\{
                                                \begin{array}{ll}
                                                  -2+\frac{2}{d_x}+\frac{2}{D(x)}+\left[\frac{(d_x-1)}{d_x}+\frac{2(d_x-1)}{d_x\vee D(x)}\right]c(x), & \hbox{if \,\,$A\geq 0$ for all $y\sim x$;} \\
                                                  -1+\frac{1}{d_x}+\frac{1}{D(x)}+\frac{2(d_x-1)}{d_x\vee D(x)}c(x), & \hbox{if \,\,$A<0\leq B$ for all $y\sim x$;} \\
                                                  \frac{d_x-1}{d_x\vee D(x)}c(x), & \hbox{if\,\,$B<0$ for all $y\sim x$.}
                                                \end{array}
                                              \right.
\end{equation*}
\end{Rmk}

\subsection{Weighted graphs}The preceding considerations readily
extend to  weighted graphs.
\begin{Thm}
 On a weighted locally finite graph $G=(V, E)$, we have
\begin{equation}
 \kappa(x, y)\geq -2\left(1-\frac{w_{xy}}{d_x}-\frac{w_{xy}}{d_y}\right)_+.
\end{equation}
Moreover, weighted trees attain this lower bound.
\end{Thm}
\begin{Thm}
On a weighted locally finite graph $G=(V, E)$, we have
\begin{align*}
\kappa(x, y)\geq&-\left(1-\frac{w_{xy}}{d_x}-\frac{w_{xy}}{d_y}-\sum_{x_1, x_1\sim x, x_1\sim y}\frac{w_{x_1x}}{d_x}\vee\frac{w_{x_1y}}{d_y}\right)_+\\
&-\left(1-\frac{w_{xy}}{d_x}-\frac{w_{xy}}{d_y}-\sum_{x_1, x_1\sim x, x_1\sim y}\frac{w_{x_1x}}{d_x}\wedge\frac{w_{x_1y}}{d_y}\right)_++\sum_{x_1, x_1\sim x, x_1\sim y}\frac{w_{x_1x}}{d_x}\wedge\frac{w_{x_1y}}{d_y}.
\end{align*}
The inequality is sharp.
\end{Thm}
\begin{Rmk}
Notice that the term replacing the number of triangles here satisfies
$$\sum_{x_1, x_1\sim x, x_1\sim y}\frac{w_{x_1x}}{d_x}\wedge\frac{w_{x_1y}}{d_y}=m_x\wedge m_y(G).$$
\end{Rmk}
\verb"Proof:" Similar to the proof of Theorem \ref{8}, we need to understand the following two terms,
\begin{align*}
 A_w:=&1-\frac{w_{xy}}{d_x}-\frac{w_{xy}}{d_y}-\sum_{x_1, x_1\sim x, x_1\sim y}\frac{w_{x_1x}}{d_x}\vee\frac{w_{x_1y}}{d_y},\\
B_w:=&1-\frac{w_{xy}}{d_x}-\frac{w_{xy}}{d_y}-\sum_{x_1, x_1\sim x, x_1\sim y}\frac{w_{x_1x}}{d_x}\wedge\frac{w_{x_1y}}{d_y}.
\end{align*}
Only the transfer plan in the case $A_w<0\leq B_w$ needs a more careful discussion. $\hfill\Box$
\begin{Thm}
 On a weighted locally finite graph $G=(V, E)$, we have for any neighboring $x, y$,
\begin{equation}
 \kappa(x, y)\leq \sum_{x_1, x_1\sim x, x_1\sim y}\frac{w_{x_1x}}{d_x}\wedge\frac{w_{x_1y}}{d_y}.
\end{equation}
\end{Thm}

\section{Curvature dimension inequalities}\label{91}

In this section, we establish curvature dimension inequalities on
locally finite graphs. A very interesting one is the inequality under
the condition $\kappa\geq k>0$. Curvature dimension inequalities on
locally finite graphs are studied in Lin-Yau \cite{LinYau}. We first
state a detailed version of their results. Let's denote
$D_w(x):=\max_{y, y\sim x}\frac{d_y}{w_{yx}}$. Notice that on an unweighted graph, this is the $D(x)$ we used in Section \ref{84}.

\begin{Thm}\label{92}
On a weighted locally finite graph $G=(V, E)$, the Laplace operator $\Delta$ satisfies
\begin{equation}\label{88}
\Gamma_2(f, f)(x)\geq\frac{1}{2}(\Delta f(x))^2+\left(\frac{2}{D_w(x)}-1\right)\Gamma(f,f)(x).
\end{equation}
\end{Thm}
\begin{Rmk}
Since in this case we attach the weighted version of measure (\ref{85}), we get
$$Hf(x)=\frac{1}{4}\frac{1}{d_x}\sum_{y, y\sim x} \frac{w_{xy}}{d_y}\sum_{z, z\sim y}w_{yz}(f(x)-2f(y)+f(z))^2.$$
We only need to choose special $z=x$ in the second sum and then (\ref{87}) and (\ref{86}) imply the theorem.
\end{Rmk}

\subsection{Unweighted graphs} We again restrict ourselves to unweighted graphs.

We observe that the existence of triangles causes cancellations in calculating the term $Hf(x)$. This gives

\begin{Thm}\label{72}
On a locally finite graph $G=(V, E)$, the Laplace operator satisfies
\begin{equation}
\Gamma_2(f, f)(x)\geq \frac{1}{2}(\Delta f(x))^2+\left(\frac{1}{2}t(x)-1\right)\Gamma(f, f)(x),\label{32}
\end{equation}
where $$t(x):=\min_{y, y\sim x}\left(\frac{4}{d_y}+\frac{1}{D(x)}\sharp(x, y)\right).$$
\end{Thm}

\begin{Rmk}
Notice that if there is a vertex $y$, $y\sim x$, such that $\sharp(x, y)=0$, this will reduce to (\ref{88}).
\end{Rmk}

\verb"Proof:" Starting from (\ref{86}), the main work is to compare $Hf(x)$ with
\begin{align*}
&\Gamma(f, f)(x)=\frac{1}{2}\frac{1}{d_x}\sum_{y, y\sim x}(f(y)-f(x))^2.
\end{align*}

First we try to write out $Hf(x)$ as
\begin{equation*}
Hf(x)=\frac{1}{4}\frac{1}{d_x}\sum_{y, y\sim x}\left[\frac{4}{d_y}(f(x)-f(y))^2+\frac{1}{d_y}\sum_{z, z\sim y, z\neq x}(f(x)-2f(y)+f(z))^2\right].
\end{equation*}

If there is a vertex $x_1$ which satisfies $x_1\sim x$, $x_1\sim y$, we have
\begin{align}
&\frac{1}{d_y}(f(x)-2f(y)+f(x_1))^2+\frac{1}{d_{x_1}}(f(x)-2f(x_1)+f(y))^2\nonumber\\
\geq& \frac{1}{D(x)}[(f(x)-f(y))^2+(f(y)-f(x_1))^2+2(f(x)-f(y))(f(x_1)-f(y))\nonumber\\
&+(f(x)-f(x_1))^2+(f(y)-f(x_1))^2+2(f(y)-f(x_1))(f(x)-f(x_1))]\nonumber\\
=&\frac{1}{D(x)}[(f(x)-f(y))^2+4(f(y)-f(x_1))^2+(f(x)-f(x_1))^2].\nonumber\\
\geq&\frac{1}{D(x)}(f(x)-f(y))^2.\label{73}
\end{align}
So the existence of a triangle which includes $x$ and $y$ will give another term $$\frac{1}{D(x)}(f(y)-f(x))^2$$
to the sum in $Hf(x)$. Since this effect is symmetric w.r.t. $y$ and $x_1$, we can get
\begin{align*}
Hf(x)&\geq\frac{1}{4}\frac{1}{d_x}\sum_{y, y\sim x}\left(\frac{4}{d_y}+\frac{1}{D(x)}\sharp(x, y)\right)(f(y)-f(x))^2\\
&\geq t(x)\frac{1}{4}\frac{1}{d_x}\sum_{y, y\sim x}(f(y)-f(x))^2\\
&=t(x)\cdot\frac{1}{2}\Gamma(f, f)(x).
\end{align*}
Inserting this into (\ref{86}) completes the proof.
$\hfill\Box$

Recalling Theorem \ref{89} and the subsequent discussion, we get the following curvature dimension inequalities on graphs with positive Ollivier-Ricci curvature.

\begin{Cor}
On a locally finite graph $G=(V, E)$, if $\kappa(x, y)>0$, then we have
\begin{equation}
\Gamma_2(f, f)(x)\geq \frac{1}{2}(\Delta f(x))^2+\left(\frac{5}{2D(x)}-1\right)\Gamma(f, f)(x).
\end{equation}
\end{Cor}

\begin{Cor}\label{78}
On a locally finite graph $G=(V, E)$, if $\kappa(x, y)\geq k>0$, then we have
\begin{equation}
\Gamma_2(f, f)(x)\geq \frac{1}{2}(\Delta f(x))^2+\left(\frac{1}{2}\min_{y, y\sim x}\left\{\frac{4}{d_y}+\frac{\lceil kd_x\vee d_y\rceil}{D(x)}\right\}-1\right)\Gamma(f, f)(x).
\end{equation}
\end{Cor}
\begin{Rmk}
Observe that a rough inequality in this case is
\begin{equation*}
\Gamma_2(f, f)(x)\geq \frac{1}{2}(\Delta f(x))^2+\left(\frac{2}{D(x)}+\frac{kd_x}{2D(x)}-1\right)\Gamma(f, f)(x).
\end{equation*}
Comparing this one with (\ref{88}), we see that positive $\kappa$ increases the curvature function here.
\end{Rmk}

\begin{Rmk}
We point out that the condition $\kappa(x, y)\geq k>0$ implies that the diameter of the graph is bounded by $\frac{2}{k}$ (see Proposition 23 in Ollivier \cite{Oll}).
So in this case the graph is a finite one.
\end{Rmk}

Let us revisit the example of a complete graph ${\cal K}_n$ $(n\geq 2)$ with $n$ vertices. Recall in Example \ref{71}, we know
\begin{equation*}
\kappa(x, y)=\frac{n-2}{n-1}, \,\,\forall \,\,x, y.
\end{equation*}
For the curvature dimension inequality on ${\cal K}_n$, Theorem \ref{72} or Corollary \ref{78} using the above $\kappa$ implies
\begin{align}
\Gamma_2(f, f)&\geq\frac{1}{2}(\Delta f)^2+\left(\frac{2}{n-1}-1+\frac{1}{2}\frac{n-2}{n-1}\right)\Gamma(f, f)\nonumber\\
&=\frac{1}{2}(\Delta f)^2+\frac{4-n}{2(n-1)}\Gamma(f, f).\label{75}
\end{align}
Moreover, the curvature term in the above inequality cannot be
larger. To see this, we calculate, using the same trick as in (\ref{73}),
\begin{align*}
Hf(x)&=\frac{1}{4(n-1)^2}\sum_{y, y\sim x}\sum_{z, z\sim x}(f(x)-2f(y)+f(z))^2\\
&=\frac{n+2}{2(n-1)}\Gamma(f, f)(x)+\frac{1}{(n-1)^2}\sum_{(x_1,\, x_2)}(f(x_1)-f(x_2))^2,
\end{align*}
where $\sum_{(x_1,\, x_2)}$ means the sum over all unordered pairs of neighbors of $x$. Recalling (\ref{86}), we get
\begin{equation}\label{74}
\Gamma_2(f, f)(x)=\frac{1}{2}(\Delta f)^2(x)+\frac{4-n}{2(n-1)}\Gamma(f, f)(x)+\frac{1}{(n-1)^2}\sum_{(x_1,\, x_2)}(f(x_1)-f(x_2))^2.
\end{equation}
For any vertex $x$, we can find a particular function $\overline{f}$,
\begin{equation}\label{94}
  \overline{f}(z)=\left\{
         \begin{array}{ll}
           2, & \hbox{when $z=x$;} \\
           1, & \hbox{when $z\sim x$,}
         \end{array}
       \right.
\end{equation}
such that the last term in (\ref{74}) vanishes, and $\Gamma(\overline{f},\overline{f})\neq 0$. This means the curvature term in (\ref{75}) is optimal for dimension parameter 2.

But the curvature term $\frac{4-n}{2(n-1)}$ behaves very differently from $\kappa$. In fact as $n\rightarrow +\infty$, $$\frac{4-n}{2(n-1)}\searrow -\frac{1}{2}\,\,\, \text{whereas}\,\,\, \kappa\nearrow 1.$$

To get a curvature dimension inequality with a curvature term which behaves like $\kappa$, it seems that we should adjust the dimension parameter. In fact, we have

\begin{Pro}
 On a complete graph ${\cal K}_n$ $(n\geq 2)$ with $n$ vertices, the Laplace operator $\Delta$ satisfies for $m\in [1, +\infty]$,
\begin{equation}
 \Gamma_2(f, f)(x)\geq\frac{1}{m}(\Delta f(x))^2+\left(\frac{4-n}{2(n-1)}+\frac{m-2}{m}\right)\Gamma(f, f)(x).
\end{equation}
Moreover, for every fixed dimension parameter $m$, the curvature term is optimal.
\end{Pro}
\verb"Proof:"  We have from (\ref{74})
\begin{align*}
 \Gamma_2(f, f)(x)=&\frac{1}{m}(\Delta f)^2(x)+\frac{4-n}{2(n-1)}\Gamma(f, f)(x)\\
&+\frac{1}{(n-1)^2}\sum_{(x_1,\, x_2)}(f(x_1)-f(x_2))^2+\left(\frac{1}{2}-\frac{1}{m}\right)(\Delta f)^2.
\end{align*}

Let us denote the sum of the last two terms by $I$. Then we have
\begin{align*}
 I=&\frac{1}{(n-1)^2}\Big\{\left(\frac{1}{2}-\frac{1}{m}\right)\sum_{y, y\sim x}(f(y)-f(x))^2+\sum_{(x_1,\,x_2)}\Big[(f(x_1)-f(x))^2+(f(x_2)-f(x))^2\\&+\left(2\left(\frac{1}{2}-\frac{1}{m}\right)-2\right)(f(x_1)-f(x))(f(x_2)-f(x))\Big]\Big\}\\
=&\frac{1}{(n-1)^2}\Big[\left(\frac{1}{2}-\frac{1}{m}\right)\sum_{y, y\sim x}(f(y)-f(x))^2+\left(1-\frac{m+2}{2m}\right)(n-2)\sum_{y, y\sim x}(f(y)-f(x))^2\\
&+\sum_{(x_1,\,x_2)}\frac{m+2}{2m}(f(x_1)-f(x_2))^2\Big]\\
=&\frac{m-2}{m}\Gamma(f, f)(x)+\frac{m+2}{2m(n-1)^2}\sum_{(x_1,\,x_2)}(f(x_1)-f(x_2))^2.
\end{align*}
This finishes the proof. $\hfill\Box$

An interesting point appears when we choose the dimension parameter $m$ of ${\cal K}_n$ as $n-1$. Then we have
\begin{equation*}
\Gamma_2(f, f)\geq \frac{1}{n-1}(\Delta f)^2+\frac{1}{2}\frac{n-2}{n-1}\Gamma(f, f),
\end{equation*}
where the curvature term is exactly $\frac{1}{2}\kappa$. From the fact that ${\cal K}_n$ could be considered as the boundary of a $(n-1)$ dimensional simplex, the $m$ we choose here seems also natural.

\begin{Rmk}
We point out another similar fact here. On a locally finite graph with maximal degree $D$ and minimal degree larger than $1$, Theorem \ref{7} and Theorem \ref{92} imply that
\begin{equation}
\kappa(x, y)\geq 2\left(\frac{2}{D}-1\right),\,\,\forall\,\, x, y,
\end{equation}
and
\begin{equation}\label{93}
\Gamma_2(f, f)\geq\frac{1}{2}(\Delta f)^2+\left(\frac{2}{D}-1\right)\Gamma(f, f),
\end{equation}
respectively. It is not difficult to see that for regular trees with degree larger than $1$, the curvature term in (\ref{93}) is optimal. (Just consider the extension of the function (\ref{94}), taking values $0$ on vertices which are not $x$ and neighbors of $x$ there.) So on regular trees, the curvature term is also exactly $\frac{1}{2}\kappa$.
\end{Rmk}

\begin{Rmk} In Erd\"os-Harary-Tutte \cite{Erdos}, they define the dimension of a graph $G$ as the minimum number $n$ such that $G$ can be embedded into a $n$ dimensional Euclidean space with every edge of $G$ having length $1$. It is interesting that by their definition, the dimension of ${\cal K}_n$ is also $n-1$ and the dimension of any tree is at most $2$.
\end{Rmk}

From the above observations, it seems natural to expect stronger relations between the lower bound of $\kappa$ and the curvature term in the curvature dimension inequality if one chooses proper dimension parameters.

\subsection{Weighted graphs} We have similar results on weighted
graphs here, with similar  proofs.
\begin{Thm}
On a weighted locally finite graph $G=(V, E)$, the Laplace operator satisfies
\begin{equation}
\Gamma_2(f, f)(x)\geq \frac{1}{2}(\Delta f(x))^2+\left(\frac{1}{2}t_w(x)-1\right)\Gamma(f, f)(x),
\end{equation}
where \begin{align*}t_w(x):=&\min_{y, y\sim x}\left\{\frac{4w_{xy}}{d_y}+\sum_{x_1, x_1\sim x, x_1\sim y}\left(\frac{w_{xy}}{d_y}\wedge\frac{w_{xx_1}}{d_{x_1}}\right)\frac{w_{x_1y}}{w_{xy}}\right\}.\end{align*}
\end{Thm}
\section*{Acknowledgement}We thank Persi Diaconis for pointing out Ollivier's notion of Ricci
curvature to us.





\begin{thebibliography}{99}
\bibitem{Bak} D. Bakry, Functional inequalities for Markov semigroups, Probability measures on groups: recent directions and trends, S. G. Dani and P. Graczyk (Editors), Tata Inst. Fund. Res., Mumbai(2006), pp. 91-147.

\bibitem{BaEm} D. Bakry and M. \'Emery, Diffusions hypercontractives (French) [Hypercontractive diffusions], S\'eminaire de probabilit\'es, XIX, 1983/84, Lecture Notes in Math. 1123, J. Az\'{e}ma and M. Yor (Editors), Springer, Berlin, 1985, pp. 177-206.

\bibitem{BaEm2} D. Bakry, and M. \'Emery, Hypercontractivit\'e de semi-groupes de diffusion. (French. English summary) [Hypercontractivity for diffusion semigroups] C. R. Acad. Sci. Paris S\'er. I Math. 299 (1984) 775-778.

\bibitem{BanJ} A. Banerjee and J. Jost, On the spectrum of the normalized graph Laplacian, Linear Algebra Appl. 428 (2008) 3015-3022.

\bibitem{BJ} F. Bauer and J. Jost, Bipartite and neighborhood graphs and the spectrum of the normalized graph Laplacian, http://arxiv.org/abs/0910.3118v3, to appear in Comm. Anal. Geom..

\bibitem{BJL} F. Bauer, J. Jost, and S. Liu, Ollivier-Ricci curvature
  and the spectrum of the normalized graph Laplace operator, preprint.


\bibitem{BS} A.-I. Bonciocat and K.-T. Sturm, Mass transportation and rough curvature bounds for discrete spaces,
J. Funct. Anal. 256 (2009) 2944-2966.

\bibitem{Chung} F. R. K. Chung,  Spectral graph theory, CBMS Regional Conference Series in Mathematics 92, American Mathematical Society, Providence, RI, 1997.

\bibitem{CY} F. R. K. Chung and S. T. Yau, Logarithmic Harnack inequalities, Math. Res. Lett. 3 (1996), 793-812.

\bibitem{DS} P. Diaconis and D. Stroock, Geometric bounds for eigenvalues of Markov chains,
Ann. Appl. Probab. 1 (1991) 36-61.

\bibitem{DSC} P. Diaconis and L. Saloff-Coste, Logarithmic Sobolev inequalities for finite Markov chains,
Ann. Appl. Probab. 6 (1996) 695-750.

\bibitem{D1} P. Diaconis, From shuffling cards to walking around the building: an introduction to modern Markov chain theory, Proceedings of the International Congress of Mathematicians, Vol. I (Berlin, 1998). Doc. Math. 1998, Extra Vol. I, 187-204.

\bibitem{D2} P. Diaconis, The Markov chain Monte Carlo revolution, Bull. Amer. Math. Soc. (N.S.) 46 (2009) 179-205.

\bibitem{DK} J. Dodziuk and L. Karp, Spectral and function theory for combinatorial Laplacians, Contemp. Math., 73, Amer. Math. Soc., Providence, RI, 1988.

\bibitem{Erdos} P. Erd\"os, F. Harary and  W. T. Tutte, On the dimension of a graph, Mathematika 12 (1965) 118-122.

\bibitem{Evans} L. C. Evans, Partial differential equations and Monge-Kantorovich mass transfer, Current developments in mathematics, 1997 (Cambridge, MA), R. Bott, A. Jaffe, D. Jerison, G. Lusztig and S. T. Yau (Editors), International Press, Boston, MA, 1999, pp.65-126.

\bibitem{Forman} R. Forman, Bochner's method for cell complexes and combinatorial Ricci curvature,
Discrete Comput. Geom. 29 (2003) 323-374.

\bibitem{Jost} J. Jost, Riemannian geometry and geometric analysis, 6th
  edition, Springer, 2011.

\bibitem{Lin} Y. Lin, Ricci curvature on graphs, 2010 John H. Barrett Memorial Lecture, http://www.math.utk.edu/barrett/2010/talks/YongLinRicci.pdf.

\bibitem{LLY} Y. Lin, L. Lu and S. T. Yau, Ricci curvature of graphs, Tohoku Math. J. 63 (2011) 605-627.

\bibitem{LinYau} Y. Lin and S. T. Yau, Ricci curvature and eigenvalue estimate on locally finite graphs, Math. Res. Lett. 17 (2010) 343-356.

\bibitem{LV} J. Lott and C. Villani, Ricci curvature for metric measure spaces via optimal transport,  Ann. of Math. 169  (2009) 903-991.

\bibitem{Ohta} S.-I. Ohta, On the measure contraction property of metric measure spaces, Comment. Math. Helv. 82 (2007) 805-828.

\bibitem{Oll} Y. Ollivier, Ricci curvature of Markov chains on metric spaces,  J. Funct. Anal. 256  (2009) 810-864.

\bibitem{Oll2} Y. Ollivier, A survey of Ricci curvature for metric spaces and Markov chains, Probabilistic approach to geometry, Adv. Stud. Pure Math. 57, M. Kotani, M. Hino and T. Kumagai (Editors), Math. Soc. Japan, Tokyo, 2010, pp. 343-381.

\bibitem{OllVill} Y. Ollivier and C. Villani, A curved Brunn--Minkowski inequality on the discrete hypercube, http://arxiv.org/abs/1011.4779v3.

\bibitem{Paeng} S.-H. Paeng, Volume, diameter of a graph and Ollivier's Ricci curvature, preprint.

\bibitem{Stone} D. A. Stone, A combinatorial analogue of a theorem of Myers, Illinois J. Math. 20 (1976), 12-21. Correction, Illinois J. Math. 20 (1976), 551-554.

\bibitem{Sturm1} K.-T. Sturm, On the geometry of metric measure spaces. I., Acta Math. 196 (2006) 65-131.

\bibitem{Sturm2} K.-T. Sturm, On the geometry of metric measure spaces. II., Acta Math. 196 (2006) 133-177.

\bibitem{V1} C. Villani, Topics in optimal transportation,  Graduate Studies in Mathematics, 58. American Mathematical Society, Providence, RI, 2003.

\bibitem{V2} C. Villani, Optimal transport, Old and new, Grundlehren der Mathematischen Wissenschaften, 338. Springer-Verlag, Berlin, 2009.

\bibitem{WS} D. J. Watts and S. H. Strogatz, Collective dynamics of 'small-world' networks. Nature 393 (1998) 440-442.
\end{thebibliography}



\end{document}